\input epsf
\documentstyle{amsppt}
\magnification=\magstep1
%\pagewidth{5.4truein}\hcorrection{0.55in}
%\pageheight{7.5truein}\vcorrection{0.75in}
\pagewidth{6.48truein}%\hcorrection{0.55in}
\pageheight{9.0truein}%\vcorrection{0.75in}
\TagsOnRight
\NoRunningHeads
\catcode`\@=11
\def\logo@{}
\footline={\ifnum\pageno>1 \hfil\folio\hfil\else\hfil\fi}
\topmatter
\title Plane partitions II: ${\bold 5\frac{\bold 1}{\bold 2}}$ 
symmetry classes
\endtitle
\author Mihai Ciucu {\rm and} Christian Krattenthaler\endauthor

\affil
  Institute for Advanced Study\\
  School of Mathematics\\
  Princeton, NJ 08540, USA\\
\\
and\\
\\
Institut f\"ur Mathematik der Universit\"at Wien\\
Strudlhofgasse 4, A-1090 Wien, Austria
\endaffil
\date July, 1998\enddate
\abstract
We present new, simple proofs for the enumeration of five of the ten symmetry classes
of plane partitions contained in a given box. Four of them are derived from a simple
determinant evaluation, using combinatorial arguments. The previous proofs of these
four cases were quite complicated.
For one more symmetry class we give
an elementary proof in the case when two of the sides of the box are equal.
Our results include simple evaluations of the determinants 
$\det\left(\delta_{ij}+{x+i+j\choose i}\right)_{0\leq i,j\leq n-1}$ and 
$\det\left({x+i+j\choose 2j-i}\right)_{0\leq i,j\leq n-1}$, notorious in plane
partition enumeration, whose previous
evaluations were quite intricate.

\endabstract
\endtopmatter
\document

\def\mysec#1{\bigskip\centerline{\bf #1}\message{ * }\nopagebreak\par}

\def\myref#1{\item"{[{\bf #1}]}"} 
 
\def\pf{{\it Proof.\ }} 

\def\cite#1{\relaxnext@
  \def\nextiii@##1,##2\end@{[{\bf##1},\,##2]}%
  \in@,{#1}\ifin@\def\next{\nextiii@#1\end@}\else
  \def\next{[{\bf#1}]}\fi\next}
\def\proclaimheadfont@{\smc}

\def\pf{{\it Proof.\ }}

\define\Z{{\bold Z}}

\define\twoline#1#2{\line{\hfill{\smc #1}\hfill{\smc #2}\hfill}}
\define\twolinetwo#1#2{\line{{\smc #1}\hfill{\smc #2}}}
\define\twolinethree#1#2{\line{\phantom{poco}{\smc #1}\hfill{\smc #2}\phantom{poco}}}

\def\mypic#1{\epsffile{#1}}

\mysec{1. Introduction}

\medskip
A plane partition is an array of nonnegative integers with the property that all rows
and columns are weakly decreasing. By a well-known bijection (see \cite{9} or
\cite{18}), plane
partitions contained in an $a\times b$ rectangle and with entries at most $c$ can be
identified with lozenge tilings of a hexagon $H(a,b,c)$ with side-lengths $a$, $b$,
$c$, $a$, $b$, $c$ (in cyclic order) and angles of $120^\circ$ (a lozenge tiling of
a region on the triangular lattice is a tiling by unit rhombi with angles of
$60^\circ$ and $120^\circ$). 

In \cite{19} Stanley describes ten natural symmetry classes of plane partitions.
Strikingly, the number of elements in each symmetry class is given by a simple
product formula. The available proofs, however, are in many cases quite intricate (see 
\cite{19}, \cite{3}, \cite{13} and \cite{21}). In this paper we present simple
proofs for five symmetry classes, and for one more we give an elementary proof in the 
case when two of the numbers $a$, $b$ and $c$ are equal.

Our proofs employ Kuperberg's observation \cite{13} that the bijection mentioned in
the first paragraph maps symmetry classes of plane partitions to symmetry classes of
tilings of $H(a,b,c)$. The three basic symmetries, in the context of tilings $T$, are:

\smallpagebreak
\flushpar
$(1)$ the reflection $t:\,T\mapsto T^t$ (called {\it transposition}) in the diagonal 
joining the two vertices of
$H(a,b,c)$ where sides of lengths $a$ and $b$ meet (this assumes $a=b$),

\smallpagebreak
\flushpar
$(2)$ the {\it rotation} $r:\,T\mapsto T^r$ by $120^\circ$ around the center of $H(a,b,c)$ 
(assuming $a=b=c$), 
and

\smallpagebreak
\flushpar
$(3)$ the rotation $k:\,T\mapsto T^k$ by $180^\circ$ (called {\it
complementation}) around the center of 
$H(a,b,c)$.
 
\smallpagebreak
If a tiling is invariant under one
of these symmetries, it is called symmetric, cyclically-symmetric or
self-complementary, respectively. 

We employ simple combinatorial arguments to deduce four difficult symmetry 
classes from a determinant evaluation due to
Andrews and Burge \cite{4}, which was later generalized by Krattenthaler \cite{12},
and then proved in a very simple way by Amdeberhan \cite{1}. The main tool used in
our proofs is the Factorization Theorem for perfect matchings presented in \cite{6}.

The first of this group of four symmetry classes is the case of cyclically
symmetric plane partitions (i.e., $T^r=T$), first proved by Andrews \cite{2} (for
another proof and a $q$-version, see \cite{15}). In
fact, Andrews' result \cite{2,Theorem 8} is a generalization of this case, and it
gives a simple product formula for

$$\det\left(\delta_{ij}+{x+i+j\choose i}\right)_{0\leq i,j\leq n-1}.\tag1.1$$
Our proof also addresses this more general result, and answers thus the problem
suggested by Mills, Robbins and Rumsey \cite{16} of finding a simple solution for the
evaluation of (1.1).

The next case we treat is that of cyclically symmetric 
transposed-complementary plane partitions (i.e., $T^r=T$ and $T^t=T^k$), first proved 
by Mills, Robbins and Rumsey \cite{16}. Again, we solve the more general problem of
evaluating the determinant

$$\det\left({x+i+j\choose 2j-i}\right)_{0\leq i,j\leq n-1}.\tag1.2$$
(It is in fact this more general result that is obtained in \cite{16}.)

The last two cases in this group of four are those of cyclically symmetric 
self-\-com\-ple\-men\-ta\-ry (i.e., $T^r=T^k=T$) and totally symmetric
self-\-com\-ple\-men\-ta\-ry
(i.e., $T^t=T^r=T^k=T$) plane partitions, 
which were first proved by Kuperberg \cite{13} and  Andrews \cite{3}, respectively.

The fifth case we deal with is that of transposed-complementary plane partitions 
(i.e., $T^t=T^k$), which was first proved by Proctor \cite{17} using arguments from
representation theory. We deduce it as a simple consequence of results in
\cite{8} on the tiling generating function of certain regions on the triangular
lattice. 

Finally, the ``half'' case --- which we provide a simple proof for, based on the
aforementioned results of \cite{8}, in case two of the numbers $a$, $b$ and $c$ are
equal --- is that of self-complementary plane partitions (i.e., $T^k=T$), which was 
first proved by Stanley \cite{19} using the theory of symmetric functions.

In fact, one more ``half-case'' could be added to the ones mentioned above: if two of
the numbers $a$, $b$ and $c$ are equal, the base case (i.e., no symmetry requirements)
follows directly by specializing $k=0$ in \cite{8,Theorem 1.1(a)}.

\mysec{2. A determinant with two tiling interpretations}

\medskip
The determinant evaluation mentioned in the Introduction from which we 
will derive the first four symmetry classes is the following.

\proclaim{Theorem 2.1 (Krattenthaler \cite{12})} Let $x$, $y$ and $n$ be nonnegative
integers with $x+y>0$, and set

$$K_n(x,y)=
\left(\frac{(x+y+i+j-1)!}{(x+2i-j)!\,(y+2j-i)!}\right)_{0\leq i,j\leq n-1}.\tag2.1$$
Then we have

$$\det(K_n(x,y))=\prod_{i=0}^{n-1}\frac{i!\,(x+y+i-1)!\,(2x+y+2i)_i\,(x+2y+2i)_i}
{(x+2i)!\,(y+2i)!},\tag2.2$$
where $(a)_k:=a(a+1)\cdots(a+k-1)$ is the shifted factorial.
\endproclaim

{\it Proof} ({\smc Amdeberhan} \cite{1}). We use the fact that for any 
matrix $A=(a_{ij})_{0\leq i,j\leq n-1}$ we have

$$\det A=\frac{(\det A_0^0)(\det A_{n-1}^{n-1})-(\det A_0^{n-1})(\det A_{n-1}^0)}
{\det A_{0,n-1}^{0,n-1}},\tag2.3$$
where $A_{i_1,\dotsc,i_k}^{j_1,\dotsc,j_k}$ is the submatrix of $A$ obtained by
removing rows indexed by $i_1,\dotsc,i_k$ and columns indexed by $j_1,\dotsc,j_k$
(see e.g. \cite{11}).

Take $A=K_n(x,y)$ in (2.3). It is readily seen that the five determinants on the
right hand side can be written in the form $\det K_m(x',y')$, with $m<n$ and suitable
$x'$ and $y'$. More precisely, we obtain
$$\align
&\det K_n(x,y)=\\
&\ \frac{\det K_{n-1}(x+1,y+1)\cdot\det K_{n-1}(x,y)-
\det K_{n-1}(x+2,y-1)\cdot\det K_{n-1}(x-1,y+2)}
{\det K_{n-2}(x+1,y+1)}.
\endalign$$
It is easy to check that the expression on the right hand side of (2.2) also
satisfies the above recurrence. Thus, (2.2) follows by induction on $n$. $\square$

\medskip
We define a {\it region} to be any subset of the plane that can be obtained as the
union of finitely many unit triangles of the regular triangular lattice. In a lozenge
tiling of a region, we allow the tile positions to be weighted. The weight of a
tiling is the product of the weights of the positions occupied by lozenges. The tiling
generating function $L(R)$ of a region $R$ is the sum of the weights of all its
tilings. 

We now introduce two types of regions whose tiling generating functions turn out to be 
very closely related to the determinant in the statement of Theorem 2.1.

Let $n$ and $x$ be nonnegative integers. Consider the pentagonal region illustrated
by Figure 2.1, where the top side has length $x$, the southeastern side has length
$n$, and the western and northeastern sides follow zig-zag paths of length $2n$.
Weight the $n$ tile positions fitting in the indentations of the northeastern boundary
by $1/2$ (we indicate weightings by 1/2 in our figures by placing shaded ovals in 
the corresponding tile positions; see Figure 2.1); weight all the others by 1. Denote 
the resulting region by $A_{n,x}$.

Let $B_{n,x}$ be the region with the same boundary as $A_{n,x}$, but having the $n-1$
tile positions fitting in the indentations of the western boundary weighted by $1/2$,
and all other tile positions weighted by 1 (see Figure 2.2). 

\topinsert
\twoline{\mypic{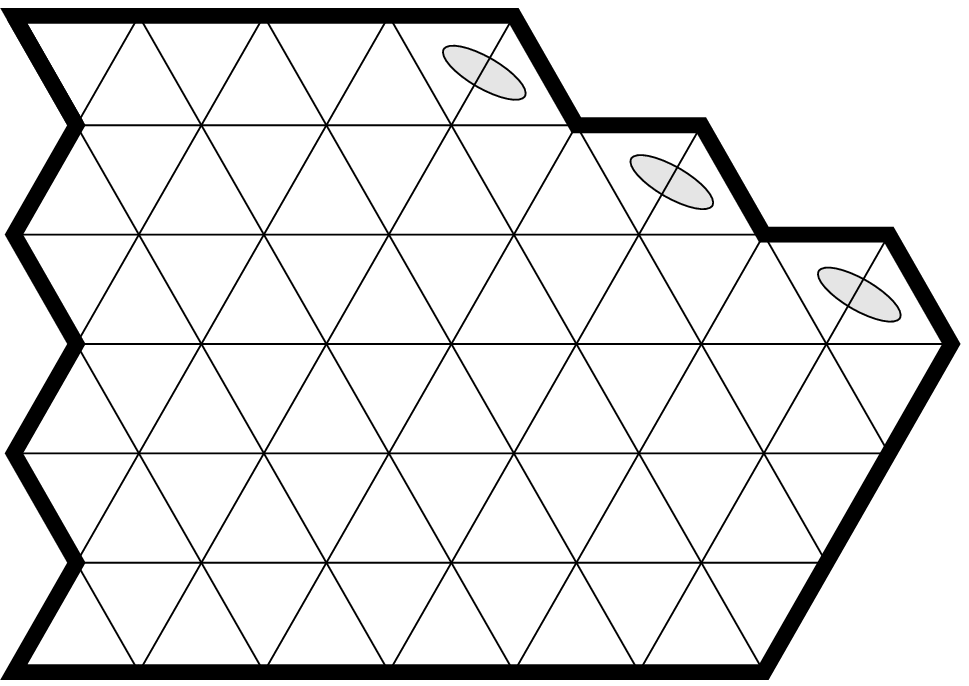}}{\mypic{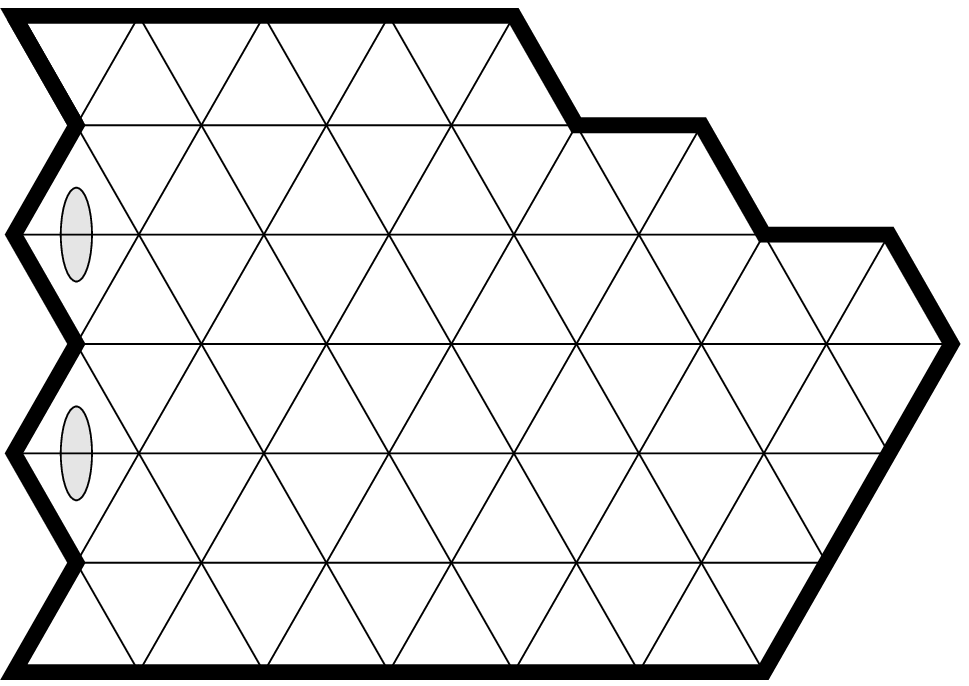}}
\twoline{Figure~2.1. {\rm $A_{3,4}$.}}{Figure~2.2. {\rm $B_{3,4}$.}}
\endinsert

The close connection between these regions and the determinant of the matrix (2.1) is 
expressed by the following result. 

\proclaim{Proposition 2.2}
$$\align
L(A_{n,x})&=\frac{1}{2^n}\det(K_n(x,0))\prod_{i=0}^{n-1}(x+3i)\tag2.4\\
L(B_{n,x})&=\frac{1}{2^n}\det(K_n(x,0))\prod_{i=0}^{n-1}(2x+3i).\tag2.5
\endalign$$
\endproclaim

\pf We use the well-known procedure of encoding tilings of a region as families of
non-intersecting lattice paths (see e.g. \cite{8, \S 4}). By this, every
tiling $T$ of $A_{n,x}$ is identified with an $n$-tuple of paths of rhombi of $T$,
each going from the western boundary of $A_{n,x}$ to its northeastern boundary. 
It follows that 
$L(A_{n,x})$ is equal to the generating function of $n$-tuples of non-intersecting
lattice paths on the square lattice taking steps north and east, starting at the
points $u_i=(i,2n-2i-1)$ and ending at $v_i=(x+2i,2n-i-1)$, $i=0,\dotsc,n-1$, where 
paths with the last step horizontal are weighted by 1/2 (the weight of a family of 
paths is the product of the weights of its elements).

It is immediate to check that the $u_i$'s and $v_j$'s satisfy the requirements
in the hypothesis of the basic theorem of Gessel and Viennot on non-intersecting
lattice paths (see e.g. \cite{20,Theorem 1.2} or \cite{10}). We obtain that
the above generating function of non-intersecting lattice paths equals 

$$\det\left((a_{ij})_{0\leq i,j\leq n-1}\right),\tag2.6$$
where $a_{ij}$ is the generating function of lattice paths from $u_i$ to $v_j$. A
straightforward calculation yields

$$\align
a_{ij}&=\frac{1}{2}{x+i+j-1\choose 2i-j}+{x+i+j-1\choose 2i-j-1}\\
&=\frac{x+3i}{2}\frac{(x+i+j-1)!}{(x-i+2j)!\,(2i-j)!}.\tag2.7
\endalign$$
Therefore, by factoring out $(x+3i)/2$ along row $i$ of the matrix in (2.6), we obtain
(2.4).

To prove (2.5) we proceed similarly. Encoding tilings as lattice paths, we obtain that
$L(B_{n,x})$ is equal to the generating function of $n$-tuples of non-intersecting
lattice paths starting and ending at the same points as above,
but now with paths having the {\it first} step {\it vertical} weighted by 1/2. It is
easy to see that in this case we have

$$\align
a_{ij}&=\frac{1}{2}{x+i+j-1\choose 2i-j-1}+{x+i+j-1\choose 2i-j}\\
&=\frac{2x+3j}{2}\frac{(x+i+j-1)!}{(x-i+2j)!\,(2i-j)!}.\tag2.8
\endalign$$
By factoring out $(2x+3j)/2$ along column $j$ we obtain (2.5). $\square$

\mysec{3. Cyclically symmetric plane partitions}

\medskip
By (2) of the Introduction, this case amounts to enumerating $r$-invariant
tilings of $H(n,n,n)$, where $n$ is a positive integer. 

We generalize this problem as follows. Consider the hexagonal region having sides of
lengths $n$, $n+x$, $n$, $n+x$, $n$, $n+x$ (in cyclic order), where $n\geq1$ and
$x\geq0$ are integers. Let $H_{n,x}$ be the region obtained from this hexagon by
removing a triangular region of side $x$ from its center, so that its vertices point
towards the shorter edges of the hexagon (this is illustrated in
Figure 3.1 for $n=4$, $x=2$). Denote by $CS(n,x)$ the number of lozenge tilings of 
$H_{n,x}$ that are
invariant under rotation by $120^\circ$. Clearly, $CS(n,0)$ is the number of 
$r$-invariant tilings of $H(n,n,n)$.

\topinsert
\twoline{\mypic{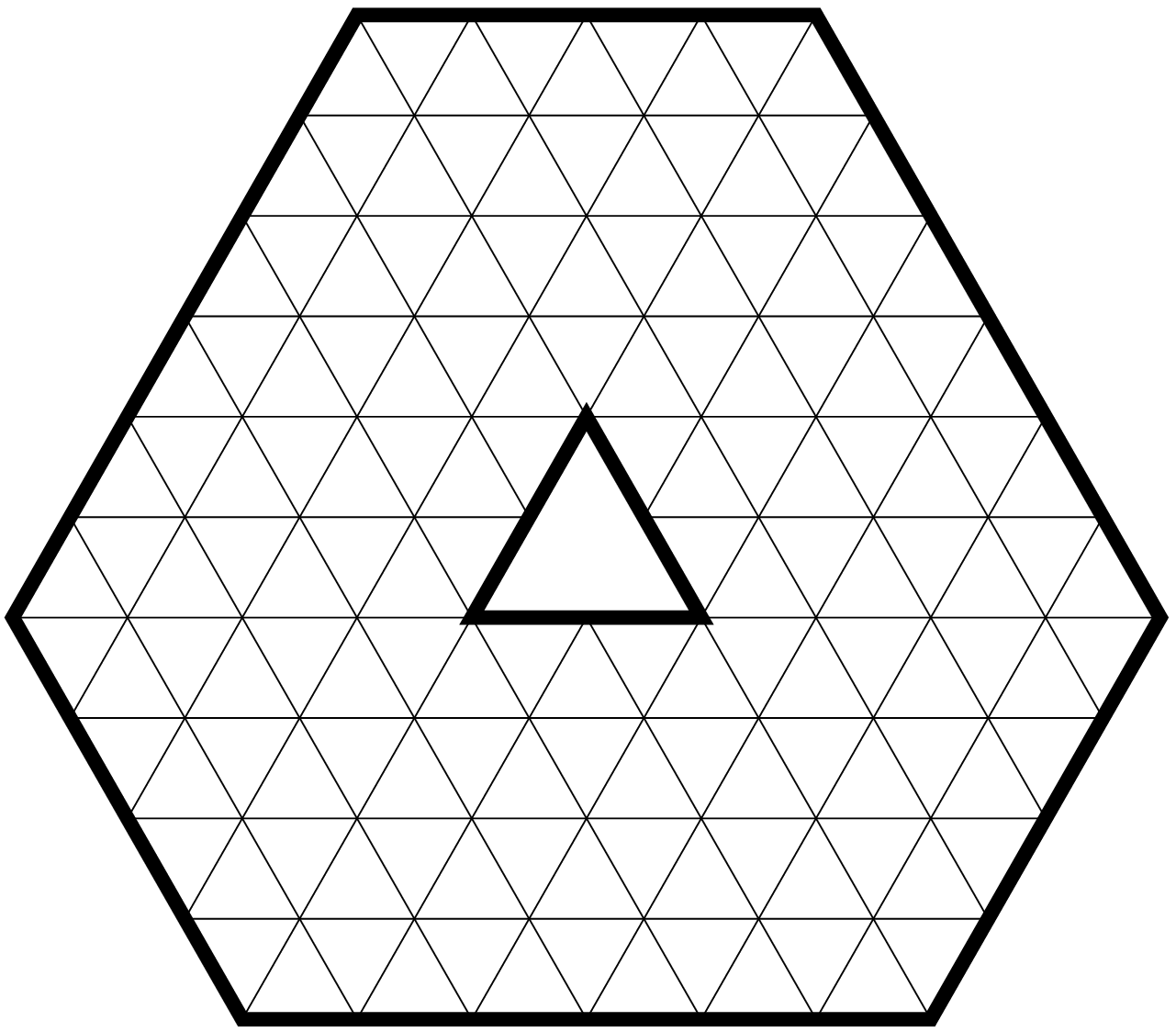}}{\mypic{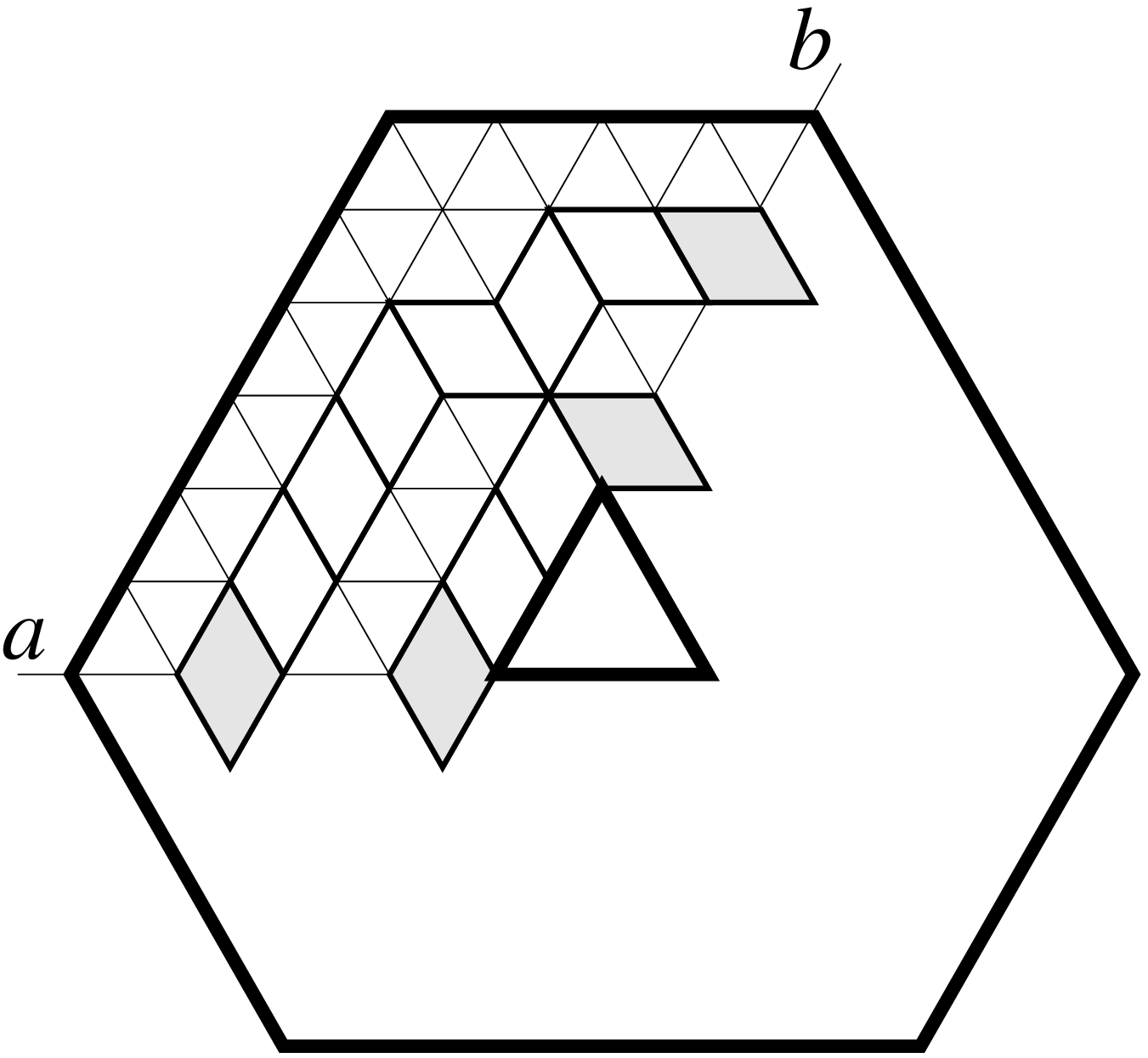}}
\twoline{Figure~3.1. {\rm \ \ \ \ \  }}{\ \ \ \ \ \ \ \ \ \ Figure~3.2. {\rm }}
\endinsert

The result below was inspired by Stembridge's proof of the special case $x=0$
(see \cite{21,Lemma 2.4}), first proved by Andrews \cite{2,Theorem 4}. 

\proclaim{Lemma 3.1} 
$$CS(n,x)=\det\left(\delta_{ij}+{x+i+j\choose i}\right)_{0\leq i,j\leq n-1}.$$
\endproclaim

\pf Let $a$ and $b$ be two lattice rays originating at two vertices of the removed 
triangle so that they determine a fundamental region $F$ for the action of $r$ on
$H_{n,x}$ (see Figure 3.2). Both $a$ and $b$ dissect $n$ tile positions in $H_{n,x}$.
Label these positions, starting with the ones closest to the removed triangle, by
$0,1,\dotsc,n-1$. 

Clearly, for any $r$-invariant tiling of $H_{n,x}$, the sets of tiles crossed
by $a$ and $b$ have the same labels. We claim that the number of tilings for which
this set of labels is  
$0\leq i_1<\cdots<i_k\leq n-1$, $0\leq k\leq n-1$, is equal to the principal minor of
the matrix $B=\left({x+i+j\choose i}\right)_{0\leq i,j\leq n-1}$ corresponding to
these labels.

Indeed, $r$-invariant tilings are determined by their intersection with the
fundamental region $F$, so such tilings with
corresponding labels $i_1,\dotsc,i_k$ can be identified with tilings of the region
$F(i_1,\dotsc,i_k)$ obtained from $F$ by removing unit triangles along $a$ and $b$ in 
positions $i_1\dotsc,i_k$. 

In turn, using the standard encoding of lozenge tilings as families of non-intersecting
lattice paths, the tilings of $F(i_1,\dotsc,i_k)$ are easily seen to be in bijection
with $k$-tuples of non-intersecting lattice paths on the square lattice, taking steps
north and east, starting at $u_\mu=(n-i_\mu-1,0)$ and ending at
$v_\mu=(n-1,x+i_\mu)$, $\mu=1,\dotsc,k$. Apply the Gessel-Viennot theorem \cite{20,
Theorem 1.2}.
The determinant corresponding to (2.6) is easily seen to be in this case precisely the 
principal minor
of $B$ corresponding to row and column indices $i_1,\dotsc,i_k$, thus proving our
claim.

We obtain that $CS(n,x)$ is equal to the sum of all principal minors of $B$, i.e., to
$\det(I+B)$. $\square$ 

\proclaim{Theorem 3.2} For $n,x\geq1$ we have

$$\align
CS(2n,2x+1)=
&\frac{n!\,(x-1)!}{(2n)!}\prod_{i=0}^{n}\frac{(x+2i)_{i+1}}{(x+n+i)!}\\
&\prod_{i=0}^{n-1}\frac{[i!]^2\,[(2x+2i+2)_{i+1}]^2\,(x+i)!\,(x+2i+1)_i}{[(2i)!]^2}
\tag3.1\endalign$$

$$\align
CS(2n-1,2x+1)=&
\frac{(x-1)!\,(2x+2n)_n}{(x+n-1)!}\\
&\prod_{i=0}^{n-1}
\frac{[i!]^2\,[(2x+2i)_{i}]^2\,(x+i)!\,(x+2i)_{i+1}\,(x+2i+1)_i}{[(2i)!]^2\,(x+n+i)!}.
\tag3.2\endalign$$

\endproclaim

\pf Lozenge tilings of $H_{2n,2x+1}$ can naturally be identified with perfect
matchings of the ``dual'' graph $G$, i.e., the graph whose vertices are the unit 
triangles
of $H_{2n,2x+1}$ and whose edges connect precisely those unit triangles that share an
edge (a perfect matching of a graph is a collection of vertex-disjoint edges
collectively incident to all vertices of the graph; we will often refer to a perfect 
matching simply as a {\it matching}). Therefore, $CS(2n,2x+1)$ is the number of
matchings of $G$ invariant under the rotation $r$ by $120^\circ$ around the center of 
$G$. 

Consider the action of the group generated by $r$ on $G$, and let $\tilde{G}$ be the
orbit graph. It follows easily that the $r$-invariant matchings of $G$ can be
identified with the matchings of~$\tilde{G}$. 

As illustrated in Figure 3.3 (for $n=3$, $x=2$), the graph $\tilde{G}$ can be 
embedded in the plane so
that it admits a symmetry axis $\ell$. Moreover, it can be readily checked that the
Factorization Theorem \cite{6,Theorem 1.2} for perfect matchings can be applied to 
$\tilde{G}$. We obtain that

\topinsert
\twolinetwo{\mypic{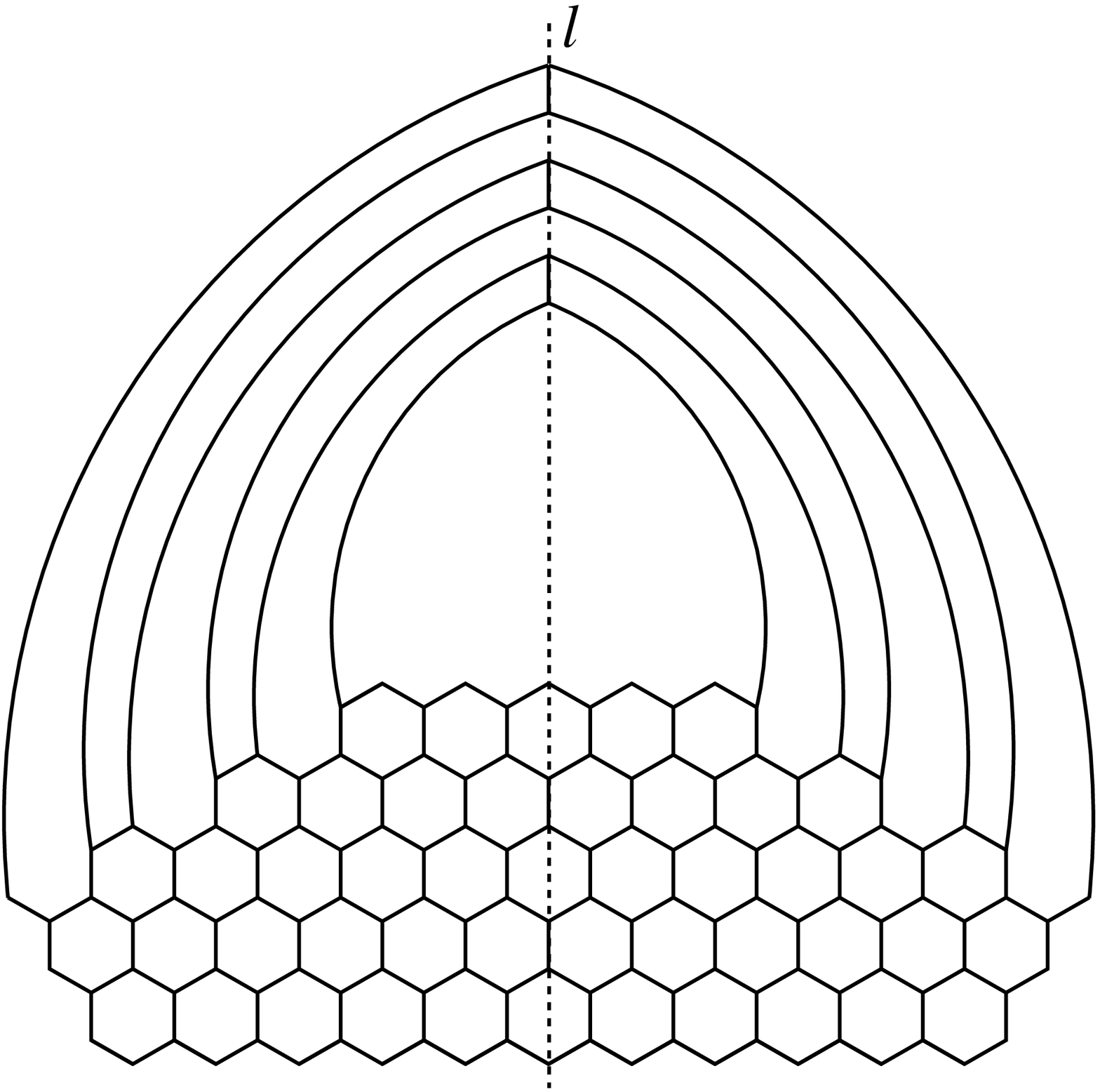}}{\mypic{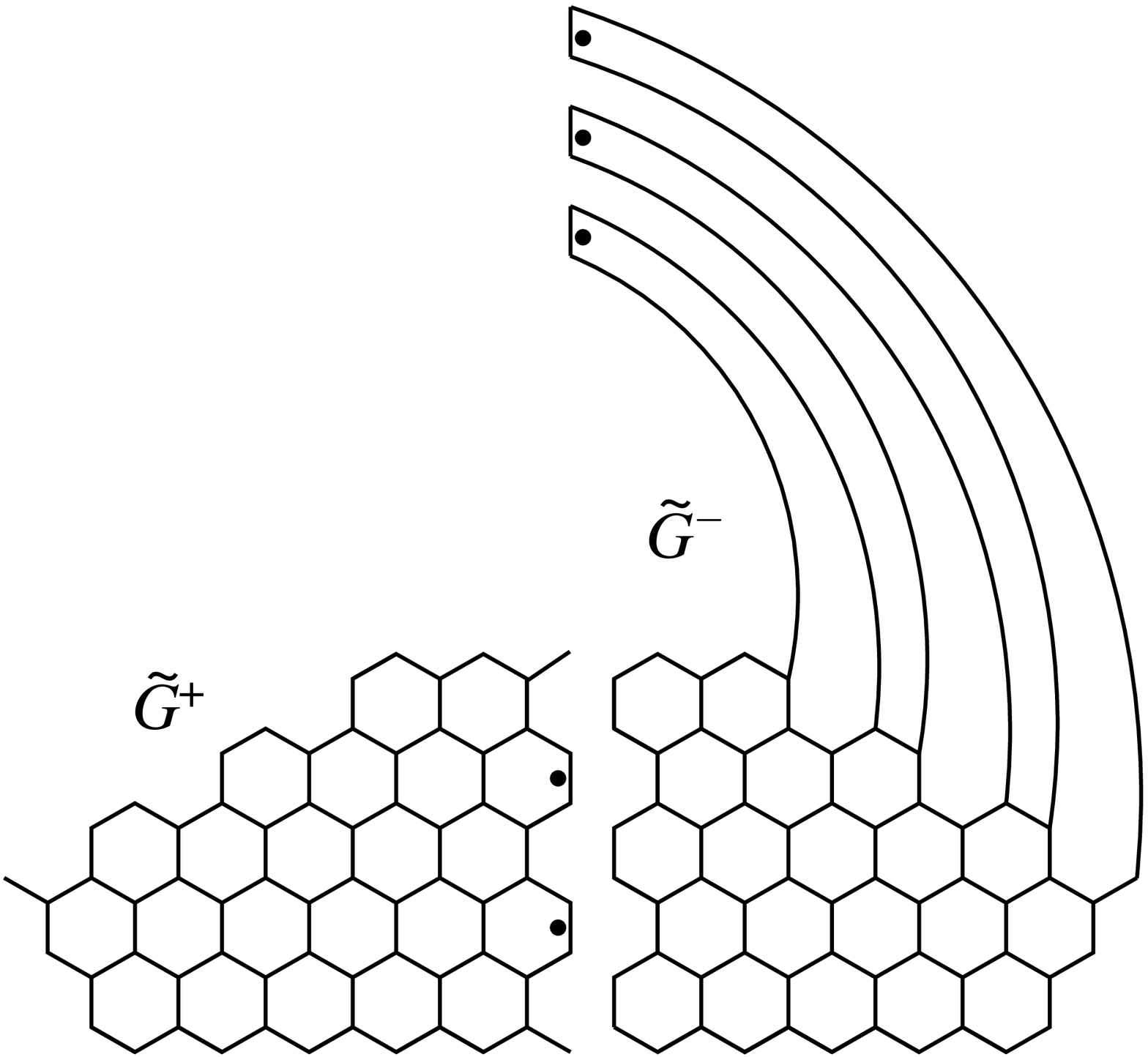}}
\twoline{Figure~3.3. {\rm \ \ \ \ \ \ \ \ \ }}
%\!\!\!\!\!\!\!\!\!\!Figure~3.3. {\rm $\tilde{G}$ for $n=3$, $x=2$.\ \ \ \ \ }}
{\ \ \ \ \ \ \ \ \ \ \ \ \ \ \ \ Figure~3.4. {\rm }}
\endinsert

\topinsert
\bigskip
\bigskip
\twolinetwo{\mypic{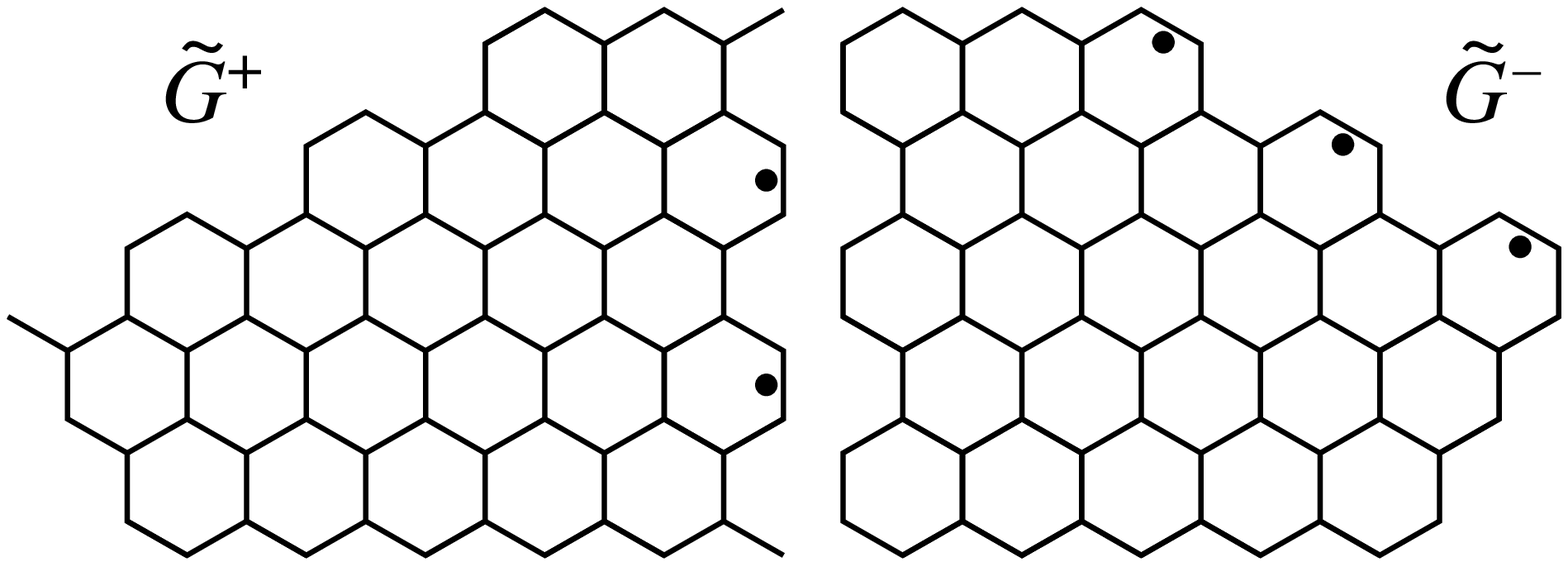}}{\mypic{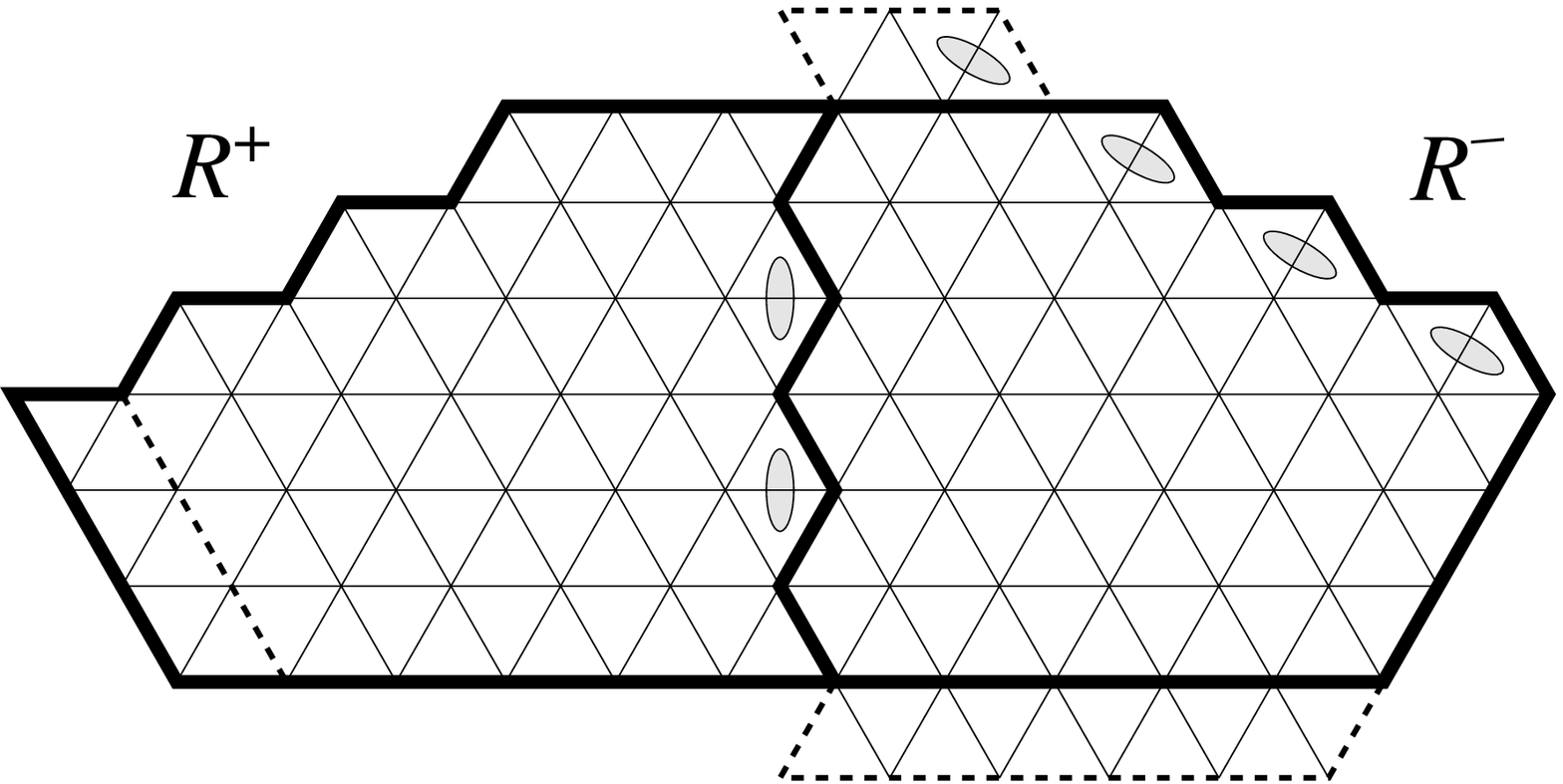}}
\twoline{Figure~3.5. {\rm \ \ \ \ \ \ \ \ \ }}
{\ \ \ \ \ \ \ \ \ \ \ \ \ \ \ \ Figure~3.6. {\rm }}
\endinsert

$$M(\tilde{G})=2^{2n}M(\tilde{G}^+)M(\tilde{G}^-),\tag3.3$$
where $M(G)$ denotes the matching generating function of $G$, and $\tilde{G}^+$ and
$\tilde{G}^-$ (illustrated in Figure 3.4) are
the connected components of the subgraph obtained from $\tilde{G}$ by deleting the 
top $2n$ edges immediately to the left of $\ell$, the bottom $2n$ edges immediately 
to the right of $\ell$, and changing the weight of the $2n-1$ 
edges along $\ell$ to 1/2 (the matching generating
function of a graph is the sum of the weights of all its perfect matchings, the
weight of a matching being the product of weights of its edges).

Clearly, the graphs $\tilde{G}^+$ and $\tilde{G}^-$ can be redrawn as shown in 
Figure 3.5. Using again the duality between matchings and tilings, we arrive at two 
regions $R^+$ and $R^-$ whose tilings can be identified, preserving weights, with the 
matchings of $\tilde{G}^+$ and $\tilde{G}^-$ (see Figure 3.6; the boundaries of $R^+$ 
and $R^-$ are shown in solid lines). However, because of 
forced tiles, it is readily seen that $L(R^+)=L(B_{n,x+1})$ and $L(R^-)=2L(A_{n+1,x})$
(compare Figure 3.6 to Figures 2.1 and 2.2; the places where the boundaries 
of $B_{n,x+1}$ and $A_{n+1,x}$ differ from those of
$R^+$ and $R^-$ are indicated by dashed lines in Figure 3.6).
Therefore, (3.3) can be rewritten as

$$CS(2n,2x+1)=2^{2n+1}L(A_{n+1,x})L(B_{n,x+1}).\tag3.4$$
By Proposition 2.2 and Theorem 2.1, the above equality yields an explicit product 
formula for $CS(2n,2x+1)$. After some manipulation one arrives at (3.1).

To prove (3.2) we proceed similarly. Take $G$ to be the graph dual to $H_{2n-1,2x+1}$,
construct the orbit graph $\tilde{G}$ as above and apply the Factorization Theorem to
$\tilde{G}$. One obtains

$$M(\tilde{G})=2^{2n-1}M(\tilde{G}^+)M(\tilde{G}^-)\tag3.5$$ 
(the change in the exponent of 2 is due the fact that the ``width'' of $\tilde{G}$ ---
cf. \cite{6}, half the number of vertices on $\ell$ --- is now $2n-1$).

\topinsert
\centerline{\mypic{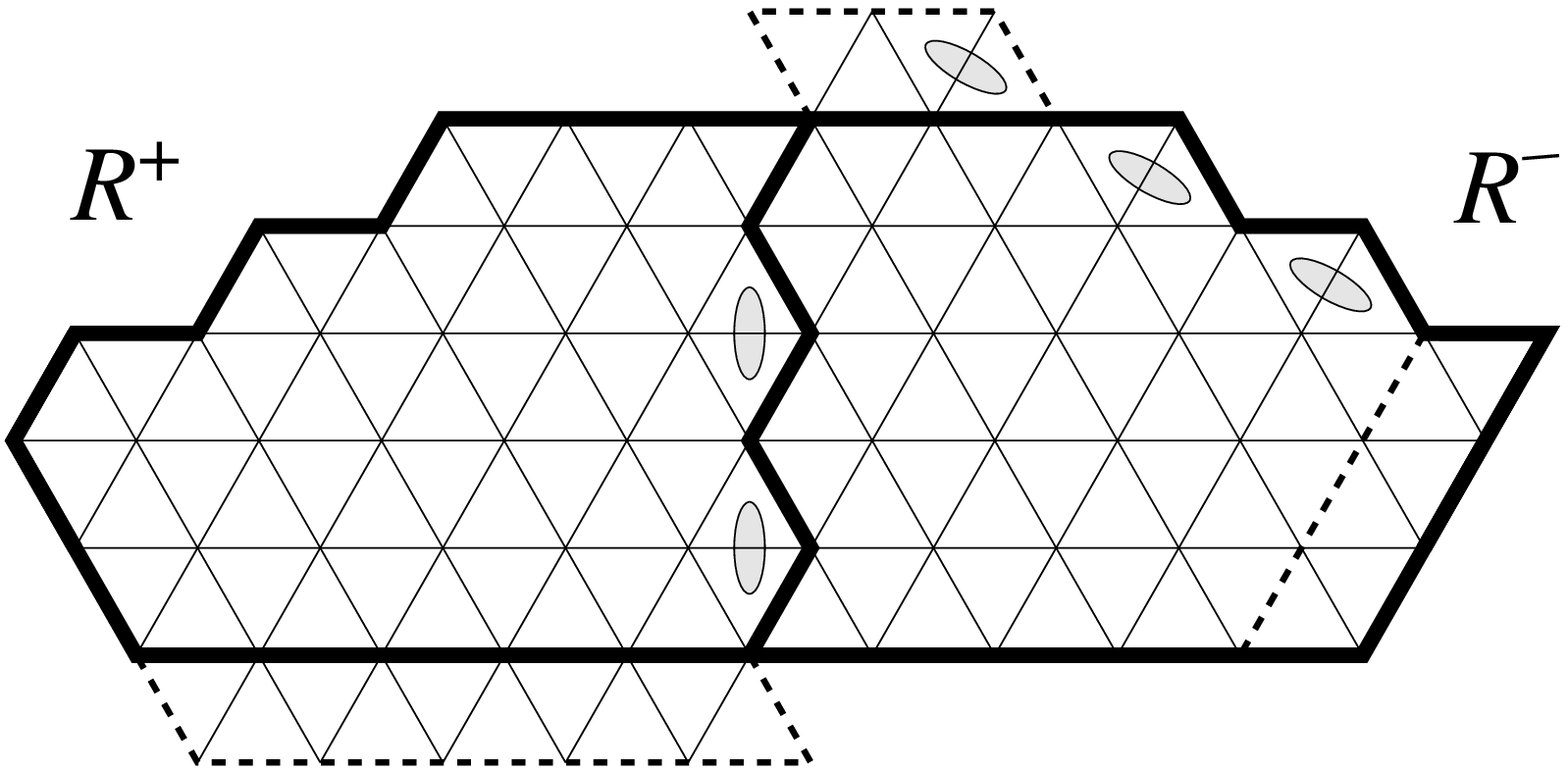}}
\centerline{Figure~3.7. {\rm }}
\endinsert

The regions $R^+$ and $R^-$ dual to $\tilde{G}^+$ and $\tilde{G}^-$ satisfy this time
$L(R^+)=L(B_{n,x+1})$ and $L(R^-)=2L(A_{n,x})$ (Figure 3.7 illustrates the case
$n=3$, $x=2$). Therefore, (3.5) implies

$$CS(2n-1,2x+1)=2^{2n}L(A_{n,x})L(B_{n,x+1}).\tag3.6$$
This provides, by Proposition 2.2 and Theorem 2.1, a product formula for
$CS(2n-1,2x+1)$, which one easily brings to the form (3.2). $\square$

\medskip
By Lemma 3.1, for fixed $n$, the expressions on the right hand side in (3.1) and (3.2)
are polynomials in $x$. Define $P_{2n}(x)$ and $P_{2n-1}(x)$ to be the polynomials on
the right hand side in (3.1) and (3.2), respectively. 

\proclaim{Corollary 3.3} With the above definition of the polynomials $P_n$, for all 
$n\geq1$ we have

$$CS(n,x)=P_{n}\left(\frac{x-1}{2}\right)\tag3.7$$
as polynomials in $x$.
\endproclaim

\pf By Theorem 3.2, (3.7) holds if $x$ is odd and $x\geq3$. Since the two sides of
(3.7) are polynomials (the left hand side by Lemma 3.1), they must be equal. $\square$

\medskip
\flushpar
{\smc Remark 3.4.} By Lemma 3.1 and Corollary 3.3 we obtain an expression for

$$\det\left(\delta_{ij}+{x+i+j\choose i}\right)_{0\leq i,j\leq n-1}$$
as a product of linear factors in $x$.
This is equivalent to Theorem 8 of \cite{2}.

\mysec{4. Cyclically symmetric transposed-complementary plane partitions}

\medskip
By (1) and (2) of the Introduction, this case is equivalent to counting tilings of
$H(n,n,n)$ that are invariant under the rotation $r$ and the reflection $t'$ across a
symmetry axis of $H(n,n,n)$ not containing any of its vertices. More
generally, we determine  the number $CSTC(n,x)$ of $r,t'$-invariant tilings of the 
regions $H_{n,x}$ (defined at the beginning of Section 3). It is easy to see that 
$H_{n,x}$ has no such tilings unless $n$ and $x$ are both even.

\topinsert
\centerline{\mypic{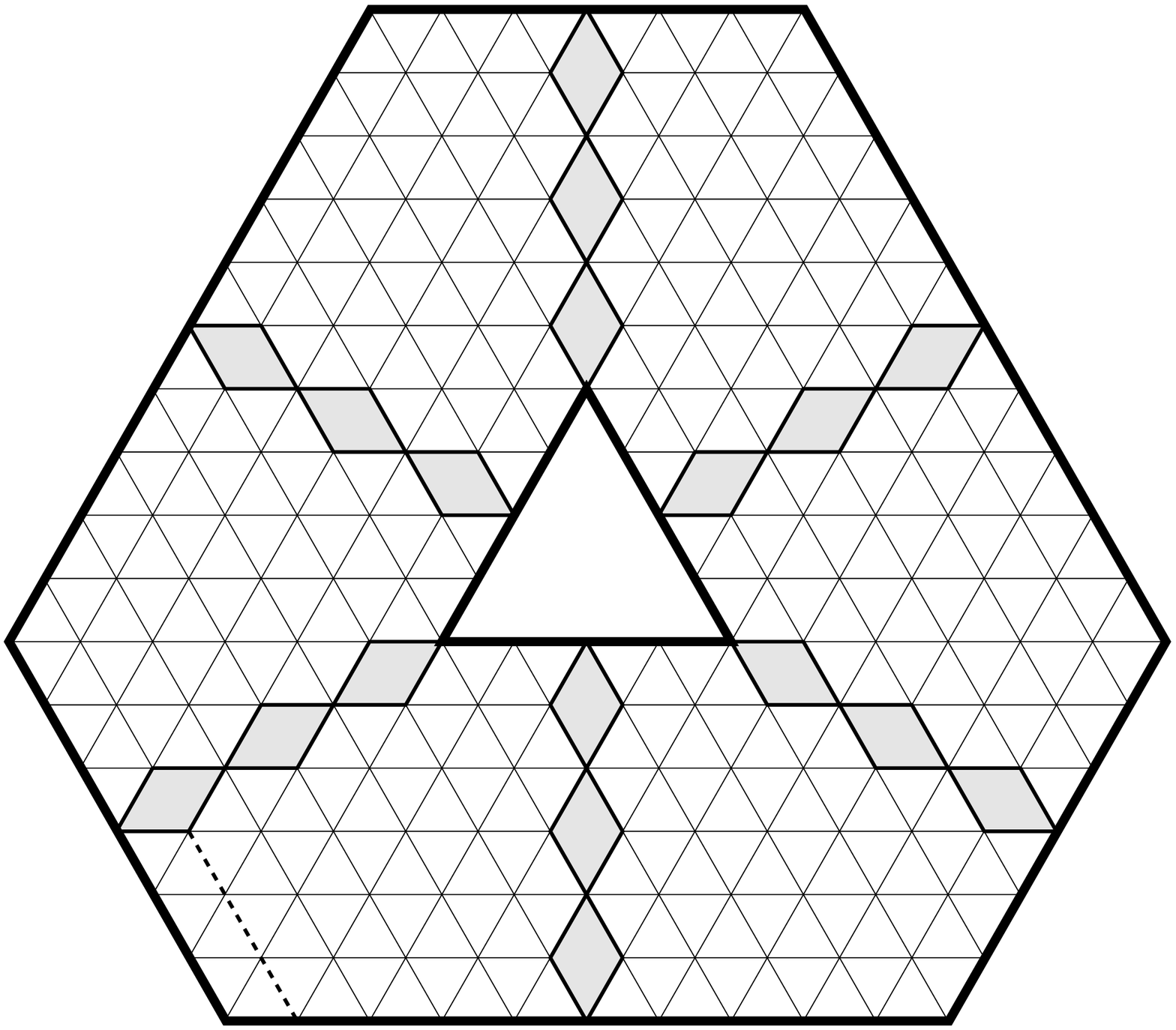}}
\centerline{Figure~4.1. {\rm }}
\endinsert

Define the region $C_{n,x}$ to be the region having the same boundary as $A_{n,x}$
(see Figure 2.1), but with all tile positions weighted by 1.

\proclaim{Lemma 4.1} $CSTC(2n,2x)=L(C_{n,x})$.
\endproclaim

\pf Suppose $T$ is an $r,t'$-invariant tiling of $H_{2n,2x}$. It follows that $T$ is 
invariant under reflection in the three symmetry axes of $H_{2n,2x}$. This implies
that in $T$ the $6n$ tile positions along these symmetry axes are occupied by
lozenges (see Figure 4.1). 
The set of these $6n$ lozenges disconnects $H_{2n,2x}$ in six congruent pieces.
Removing $n$ forced lozenges from one of these pieces one obtains a region 
congruent to $C_{n,x}$ (this is indicated by the dotted line in Figure 4.1). 
The group generated by $r$ and $t'$ acts transitively on the
set of these pieces. Therefore, the restriction of $T$ to one of the pieces gives a
bijection between $r,t'$-invariant tiling of $H_{2n,2x}$ and tilings of $C_{n,x}$.
$\square$

\proclaim{Theorem 4.2}
$$2\frac{CSTC(2n+2,2x)}{CSTC(2n,2x)}=\frac{CS(2n+1,2x)}{CS(2n,2x)}.\tag4.1$$
\endproclaim

\pf We deduce (4.1) by working out the analogs of (3.4) and (3.6) for the case when
the second argument on the left hand side is even. 

We proceed along the same lines as in the proof of Theorem 3.2. Let $G$ be the graph
dual to $H_{2n,2x}$, and let $\tilde{G}$ be the orbit graph of the action of
$\langle r\rangle$ on $G$. As in the proof of Theorem 3.2, we can embed 
$\tilde{G}$ in the plane so that it admits a symmetry axis, and we can apply the
Factorization Theorem of \cite{6}. This expresses the number of perfect matchings of
$\tilde{G}$ as a product involving the matching generating functions of two subgraphs.
These two subgraphs can be redrawn in the plane such that they are the dual graphs of
two regions $R_1^+$ and $R_1^-$ on the triangular lattice. For $n=2$, $x=1$, these 
regions are illustrated in Figure 4.2 (their boundaries are shown in solid lines). 

Therefore, since $M(\tilde{G})=CS(2n,2x)$, we can phrase the result of applying the 
Factorization Theorem to $\tilde{G}$ as

$$CS(2n,2x)=2^{2n}L(R_1^+)L(R_1^-).\tag4.2$$
However, by removing the $n$ forced lozenges along the left boundary of $R_1^+$, we
are left with a region congruent to $C_{n,x}$ (this is indicated by the dotted line
in Figure 4.2). Thus, (4.2) implies

\topinsert
\twolinetwo{\mypic{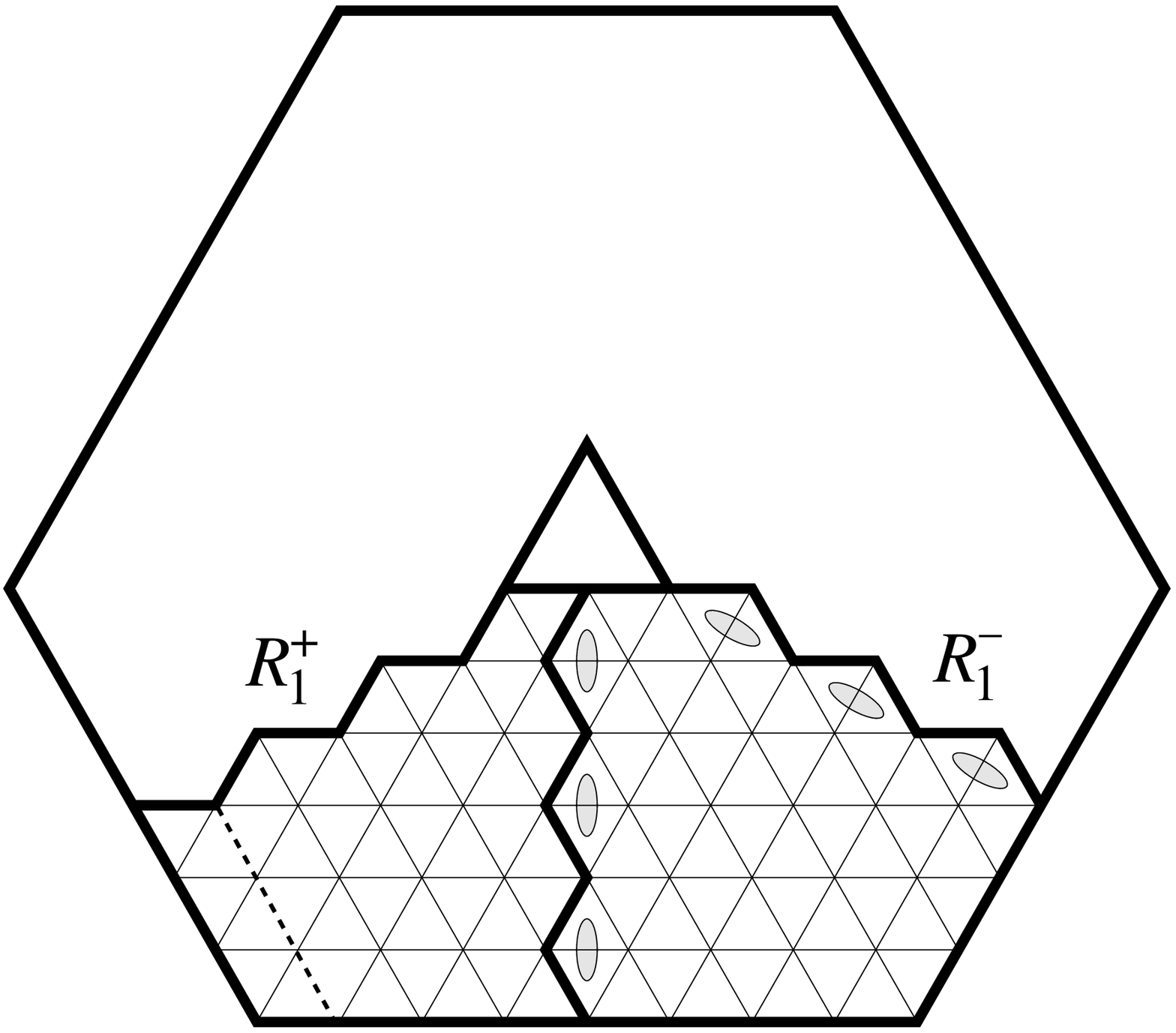}}{\mypic{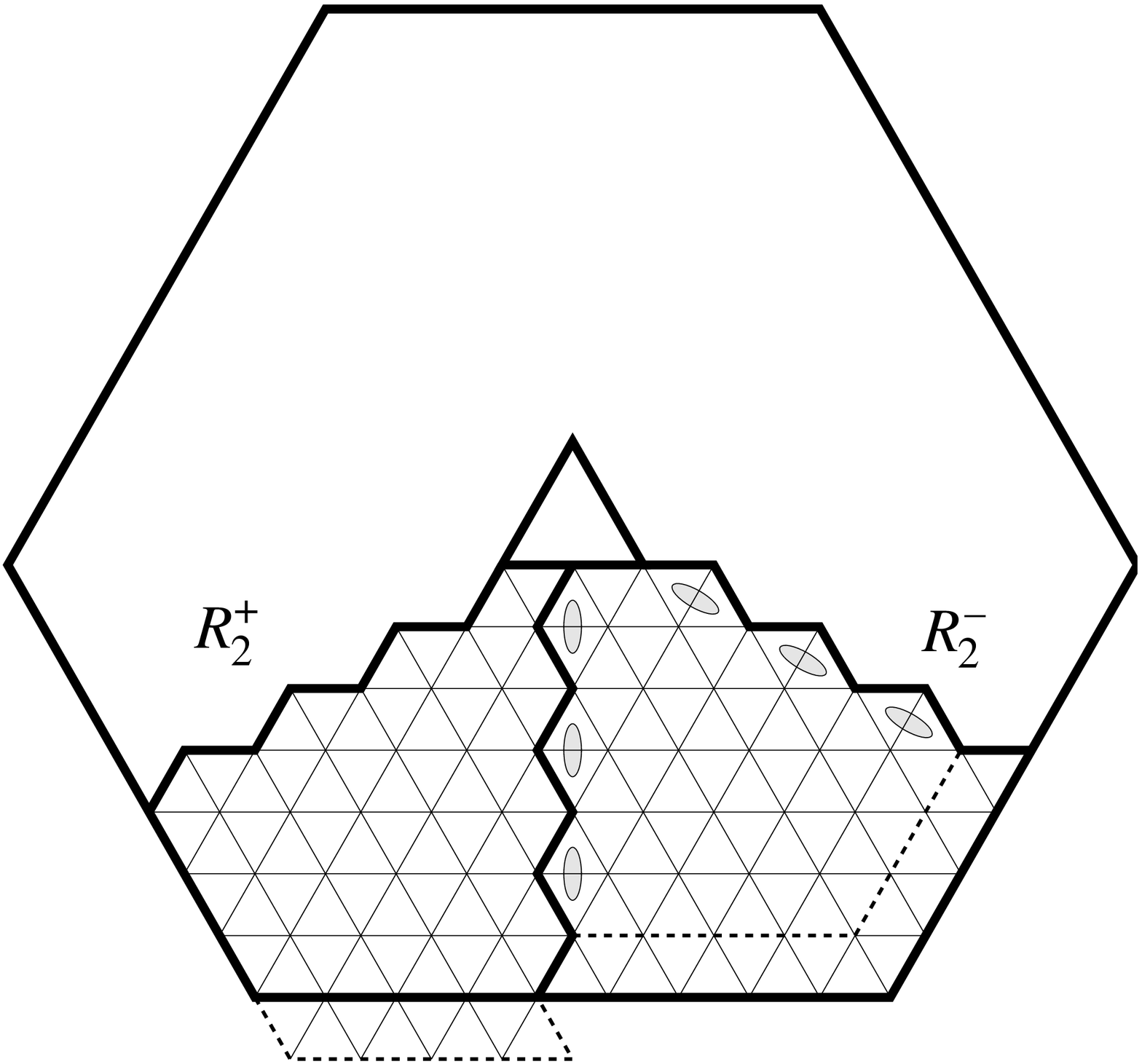}}
\twoline{Figure~4.2. {\rm \ \ \ \ \ \ \ \ \ }}
{\ \ \ \ \ \ \ \ \ \ \ \ \ \ \ \ Figure~4.3. {\rm }}
\endinsert

$$CS(2n,2x)=2^{2n}L(C_{n,x})L(R_1^-).\tag4.3$$

Similarly, starting from the graph dual to $H_{2n+1,2x}$, considering its orbit graph
under the action of $\langle r\rangle$ and applying the Factorization Theorem to it,
we obtain after rephrasing everything in terms of tilings that

$$CS(2n+1,2x)=2^{2n+1}L(R_2^+)L(R_2^-),\tag4.4.$$
for two regions $R_2^+$ and $R_2^-$ which are illustrated in Figure 4.3 for
$n=2$, $x=1$ (their boundaries are shown in solid lines). However, $R_2^+$ is
congruent to the region obtained from $C_{n+1,x}$ after removing the $n+x$ forced
lozenges along its base. Moreover, the region obtained from $R_2^-$ by removing all
forced lozenges is isomorphic to $R_1^-$. Therefore, (4.4) becomes

$$CS(2n+1,2x)=2^{2n+1}L(C_{n+1,x})L(R_1^-).\tag4.5.$$
Dividing (4.3) and (4.5) side by side and using Lemma 4.1 we obtain (4.1). $\square$

\proclaim{Corollary 4.3}
$$CSTC(2n,2x)=\frac{1}{2^n}\prod_{k=0}^{n-1}\frac{CS(2k+1,2x)}{CS(2k,2x)}.\tag4.6$$
$($By Corollary 3.3, this provides an explicit formula for $CSTC(2n,2x)$.$)$
\endproclaim

\pf Take the side by side product of (4.1) for $n=0,1,\dotsc,n-1$. $\square$

\medskip
\flushpar
{\smc Remark 4.4.} Using the standard encoding of lozenge tilings as
families of non-intersecting lattice paths, and then employing the Gessel-Viennot
theorem \cite{20, Theorem 1.2}, it is easy to see that

$$L(C_{n,x})=\det\left({x+i+j\choose 2j-i}\right)_{0\leq i,j\leq n-1}.\tag4.7$$
Therefore, by (4.7), Lemma 4.1 and Corollary 4.3 we obtain an expression for

$$\det\left({x+i+j\choose 2j-i}\right)_{0\leq i,j\leq n-1}$$
as a product of linear factors in $x$. Such a formula was first proved by Mills,
Robbins and Rumsey in \cite{16,Theorem 7}.

\medskip
\flushpar
{\smc Remark 4.5.} Following the notation of \cite{5}, let

$$Z_n(x)=\det\left(\delta_{ij}+{x+i+j\choose i}\right)_{0\leq i,j\leq n-1},$$

$$T_n(x)=\det\left({x+i+j\choose 2j-i}\right)_{0\leq i,j\leq n-1},$$

$$R_n(x)=
\det\left({x+i+j\choose 2i-j}+2{x+i+j+2\choose 2i-j+1}\right)_{0\leq i,j\leq n-1}.$$

Encode the tilings of the region  $R_1^-$ in (4.3) as $n$-tuples of non-intersecting
paths of rhombi going from the western boundary to the northeastern boundary of $R_1^-$
(see Figure 4.2). Identify, as usual, these paths of rhombi with lattice paths on
$\Z^2$. Apply the Gessel-Viennot
theorem on non-intersecting lattice paths \cite{20, Theorem 1.2}. 
It is easy to see that the
$(i,j)$-entry of the Gessel-Viennot matrix $M$ is in this case

$$M_{ij}={x+i+j\choose 2i-j}+\frac12{x+i+j\choose 2i-j-1}+\frac12{x+i+j\choose 2i-j+1}
+\frac14{x+i+j\choose 2i-j},$$
for $i,j=0,\dotsc,n-1$. A simple calculation shows that 

$$M_{ij}=\frac14\left(R_n(x)\right)_{ij}.$$

Therefore, (4.3) and (4.5) can be written as
$$\align
Z_{2n}(2x)&=T_n(x)R_n(x),\\
Z_{2n+1}(2x)&=2T_{n+1}(x)R_n(x).
\endalign$$

These are precisely relations (2.5) and (2.6) of \cite{5}, which were deduced there
from \cite{16,Theorem 5}. 

\medskip
\flushpar
{\smc Remark 4.6.} The case $x=0$ of Theorem 4.2 is the object of Theorem 6.2 of
\cite{6}.

\bigskip
\centerline{{\bf 5. Cyclically symmetric self-complementary and totally}}
\nopagebreak\centerline{{\bf symmetric self-complementary plane partitions}}

\bigskip
It is easy to see that in order for the hexagon $H(a,b,c)$ to have tilings in 
any of the two symmetry classes
mentioned in the title of this section one needs to have $a=b=c=2n$, with $n$ a
positive integer. Denote by $CSSC(2n)$ and $TSSC(2n)$ the number of tilings of
$H(2n,2n,2n)$ in the two symmetry classes, respectively.

The following result was first proved by Kuperberg \cite{13}.

\proclaim{Theorem 5.1} 
$$CSSC(2n)=\left(\prod_{i=0}^{n-1}\frac{(3i+1)!}{(n+i)!}\right)^2.\tag5.1$$
\endproclaim

\pf In \cite{7} it is shown (see \cite{7,(2.3)}) that a simple 
consequence of the Factorization Theorem for perfect matchings \cite{6,Theorem 1.2} 
is that 

$$CSSC(2n)=2^nL(A_{n,1}).\tag5.2$$
(The derivation of this result follows along the lines of the proofs of (3.4), (3.6),
(4.3) and (4.5).)
Using Proposition 2.2 and Theorem 2.1 we obtain a product formula for
$CSSC(2n)$, which is easily seen to be equivalent to (5.1). $\square$

\medskip
\flushpar
{\smc Remark 5.2.} Following the notation of \cite{5}, let

$$W_n(x)=\left({x+i+j+1\choose 2i-j+1}+
{x+i+j\choose 2i-j}\right)_{0\leq i,j\leq n-1}.$$

Based on the fact that the related determinants $Z_n(x)$ and $T_n(x)$ defined in
Remark 4.5 have close connections with plane partition enumeration problems, Andrews
and Burge suggest in \cite{5} that the same might be true for $\det\,W_n(x)$. 
Relation (5.2) allows us to give what appears to be the first such connection.

Indeed, let

$$w_n(x)=\left({x+i+j+1\choose 2i-j}+
{x+i+j\choose 2i-j-1}\right)_{0\leq i,j\leq n-1}.$$
It is readily checked that the matrix obtained from $w_n(x)$ by removing the first row
and column is precisely $W_{n-1}(x+2)$. Since the top left entry of $w_n(x)$ is 1, we
obtain that

$$\det\,W_{n-1}(x+2)=\det\,w_n(x).\tag5.3$$
A straightforward calculation reveals that $w_n(x)_{ij}=(x+3i+1)K_n(x+1,0)_{ij}$. 
Therefore, by (2.4), we deduce that

$$\det\,w_n(x)=2^nL(A_{n,x+1}).\tag5.4$$
From (5.3) and (5.4), it follows that $\det\,W_{n-1}(x+2)=2^nL(A_{n,x+1})$.
Therefore, by (5.2), we obtain that $\det\,W_{n-1}(2)=CSSC(2n)$.

\medskip
In \cite{7} there is presented a direct proof of the fact that 

$$CSSC(2n)=TSSC(2n)^2.$$
(In outline, by combinatorial arguments, an expression is derived for $CSSC(2n)$ as the
determinant of a certain matrix, which is then transformed by elementary row and
column operations to an antisymmetric matrix whose Pfaffian was previously known to
give $TSSC(2n)$).

Therefore, we obtain by Theorem 5.1 the following result, first proved by Andrews
\cite{3}.

\proclaim{Corollary 5.3} 
$$TSSC(2n)=\prod_{i=0}^{n-1}\frac{(3i+1)!}{(n+i)!}.$$
\endproclaim

\medskip
\flushpar
{\smc Remark 5.4.} By (3.4) and (3.6) we have

$$\frac{CS(2n,2x+1)}{CS(2n-1,2x+1)}=2\frac{L(A_{n+1,x})}{L(A_{n,x})}.\tag5.5$$
On the other hand, from (5.2) we deduce

$$\frac{CSSC(2n+2)}{CSSC(2n)}=2\frac{L(A_{n+1,1})}{L(A_{n,1})}.$$
This relation and (5.5) specialized to $x=1$ imply

$$\frac{CSSC(2n+2)}{CSSC(2n)}=\frac{CS(2n,3)}{CS(2n-1,3)}.\tag5.6$$
One may regard (5.6) as giving a proof of the cyclically symmetric, self-complementary
case based on the solution of the cyclically symmetric case, which was solved fifteen 
years earlier (see \cite{2} and \cite{13}).

\mysec{6. Transposed-complementary plane partitions}

\medskip
By (1) and (3) of the Introduction, this case amounts to finding the number 
$TC(a,a,2b)$ of tilings of the
hexagon $H(a,a,2b)$ that are symmetric with respect to its symmetry axis
$\ell$ perpendicular
to the sides of length $2b$ (see Figure 6.1; it is easy to see that the indicated 
form of the arguments represents the general case).  

The following result was first proved (in an equivalent form) by Proctor \cite{17}.

\proclaim{Theorem 6.1}
$$\align
&TC(a,a,2b)=\\
&  
\frac{(b+1)(b+2)^2\cdots(b+a-2)^2(b+a-1)(2b+3)(2b+5)^2\cdots(2b+2a-5)^2(2b+2a-3)}
{1\cdot2^2\cdots(a-2)^2(a-1)\cdot3\cdot5^2\cdots(2a-5)^2(2a-3)}
\endalign$$
$($the dots indicate that the bases grow by one from each factor to the next, while 
the exponents increase by one midway through and then decrease by one to the end$)$.
\endproclaim

\topinsert
\centerline{\mypic{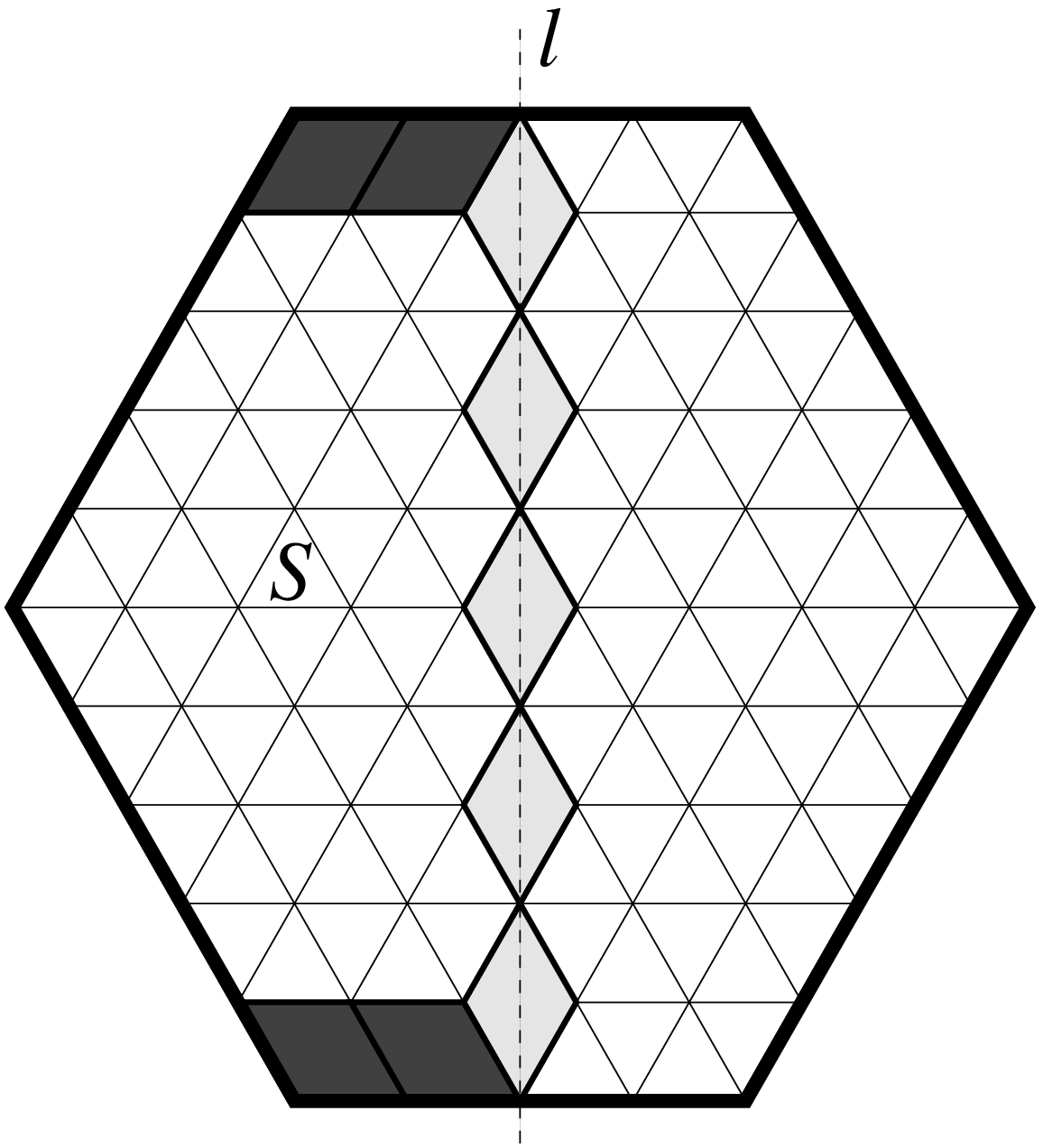}}
\centerline{Figure~6.1.}
\endinsert

\pf In any tiling $T$ of $H(a,a,2b)$ symmetric with respect to $\ell$, the $a$ tile
positions along $\ell$ are occupied by lozenges. This set of lozenges divides our
hexagon in two congruent pieces, and $T$ is determined by its restriction to the left
piece $S$, say (see Figure 6.1). Therefore, $TC(a,a,2b)$ is just the number of tilings 
of $S$. 

However, the region obtained from $S$ by removing the forced lozenges (see Figure 6.1)
is readily recognized as being a member of the family of regions 
$\bar{R}_{{\bold l},{\bold q}}(x)$ defined in \cite{8, \S 2} (here ${\bold l}$ and 
${\bold q}$ are lists of strictly increasing positive integers, and $x$ is integer). 
More precisely, $S$ is congruent to the region
$\bar{R}_{[a-1],\emptyset}(b)$, where $[n]$ denotes the list $(1,\dotsc,n)$. 

Therefore, Proposition 2.1 of \cite{8}
and formulas (1.6), (1.2) and (1.4) of \cite{8} provide an expression for $L(S)$
(hence, for  $TC(a,a,2b)$) as a
product of linear polynomials in $b$. After some manipulation, this expression becomes
the right hand side of the equality in the statement of the Theorem. $\square$

\mysec{7. Self-complementary plane partitions}

\medskip
This case amounts to enumerating tilings of $H(a,b,c)$ that are invariant under the
rotation $k$ by $180^\circ$, and it was first proved by Stanley \cite{19}. In this 
section we give a simple proof in the case when two of the numbers $a$, $b$ and $c$ 
are equal.

Let $SC(a,a,b)$ be the number of $k$-invariant tilings of $H(a,a,b)$. It is easy to
see that this number is 0 unless $a$ or $b$ is even. Define $PP(a,b,c)$ to be the
expression on the right hand side of (6.1).

\proclaim{Theorem 7.1} 
$$\align
SC(2x,2x,2y)=&PP(x,x,y)^2\tag7.1\\
SC(2x,2x,2y+1)=&PP(x,x,y)PP(x,x,y+1)\tag7.2\\
SC(2x+1,2x+1,2y)=&PP(x,x+1,y)^2.\tag7.3\\
\endalign$$
\endproclaim

\pf Following the same reasoning as in proving (3.4), (3.6), (4.3), (4.5) and (5.1), 
one sees
that the Factorization Theorem of \cite{6} can be used to express the number of
$k$-invariant tilings of $H(a,a,b)$ as a power of 2 times the tiling generating
function of a certain subregion with some tile positions weighted by 1/2 (the precise
shape of this region depends on the parities of $a$ and $b$).
Furthermore, after removing the forced lozenges from this region, the leftover
piece is readily recognized to belong to one of the families 
$R_{{\bold l},{\bold q}}(x)$ or $\bar{R}_{{\bold l},{\bold q}}(x)$ defined in 
\cite{8, \S 2}. 

More precisely, for $a=2x$, $b=2y$, we obtain that

$$L(H(2x,2x,2y))=2^xL(\bar{R}_{[x-1],[x]}(y))\tag7.4$$
(see Figure 7.1; as usual, the dotted lines indicate removal of forced lozenges). 
Similarly, we deduce 

\topinsert
\centerline{\mypic{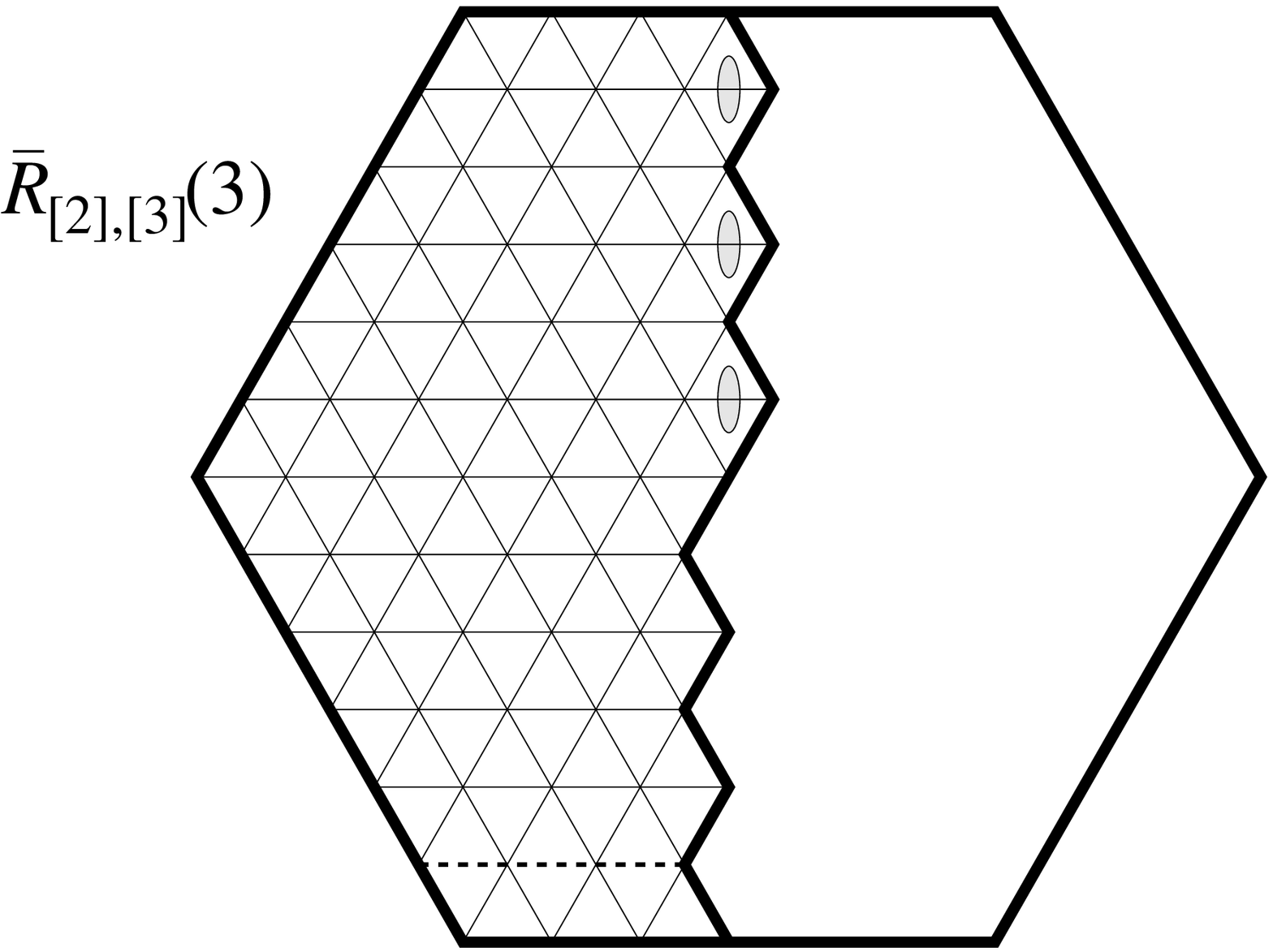}}
\centerline{Figure~7.1. $a=6$, $b=6$.}
\endinsert

\topinsert
\twoline{\mypic{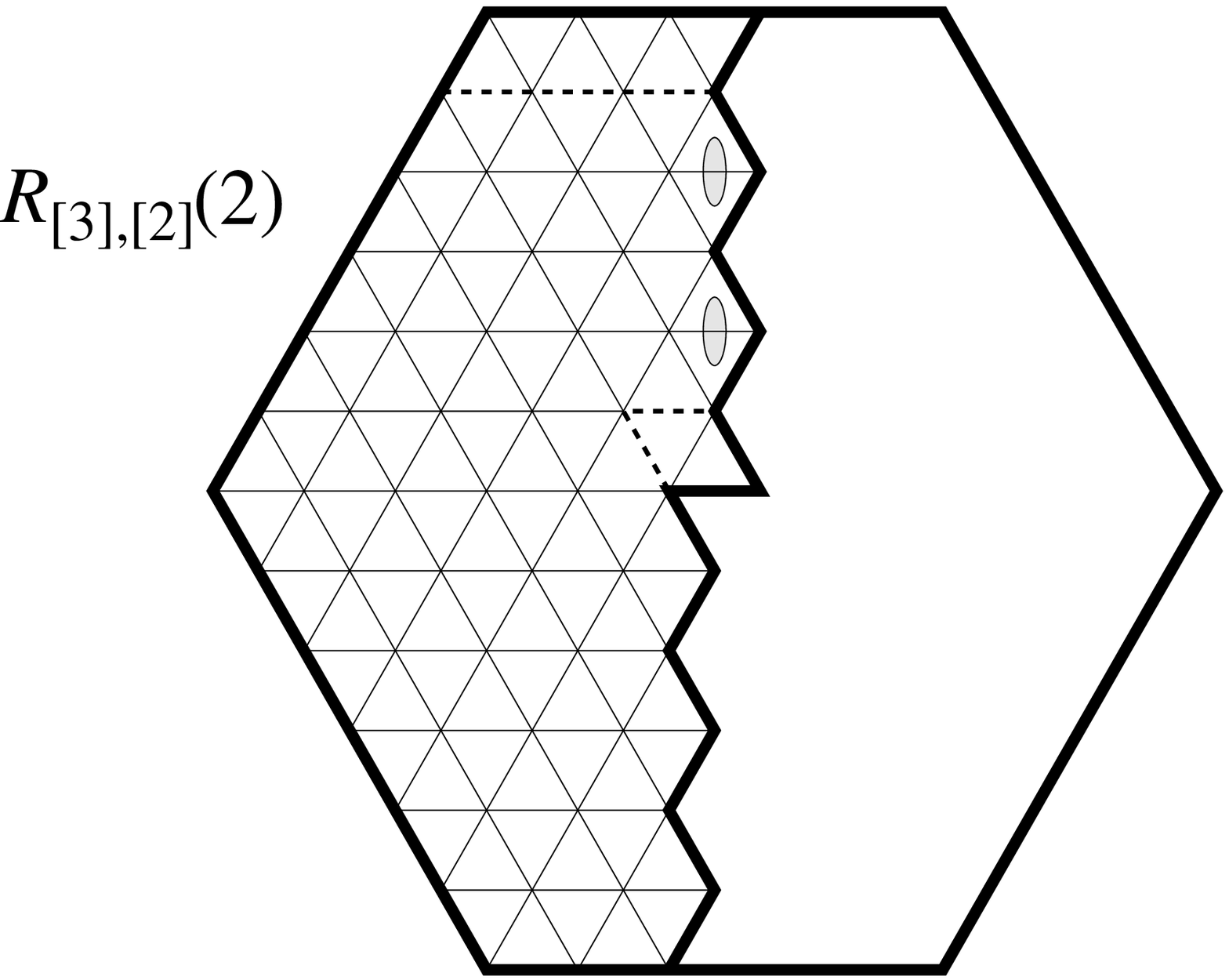}}{\mypic{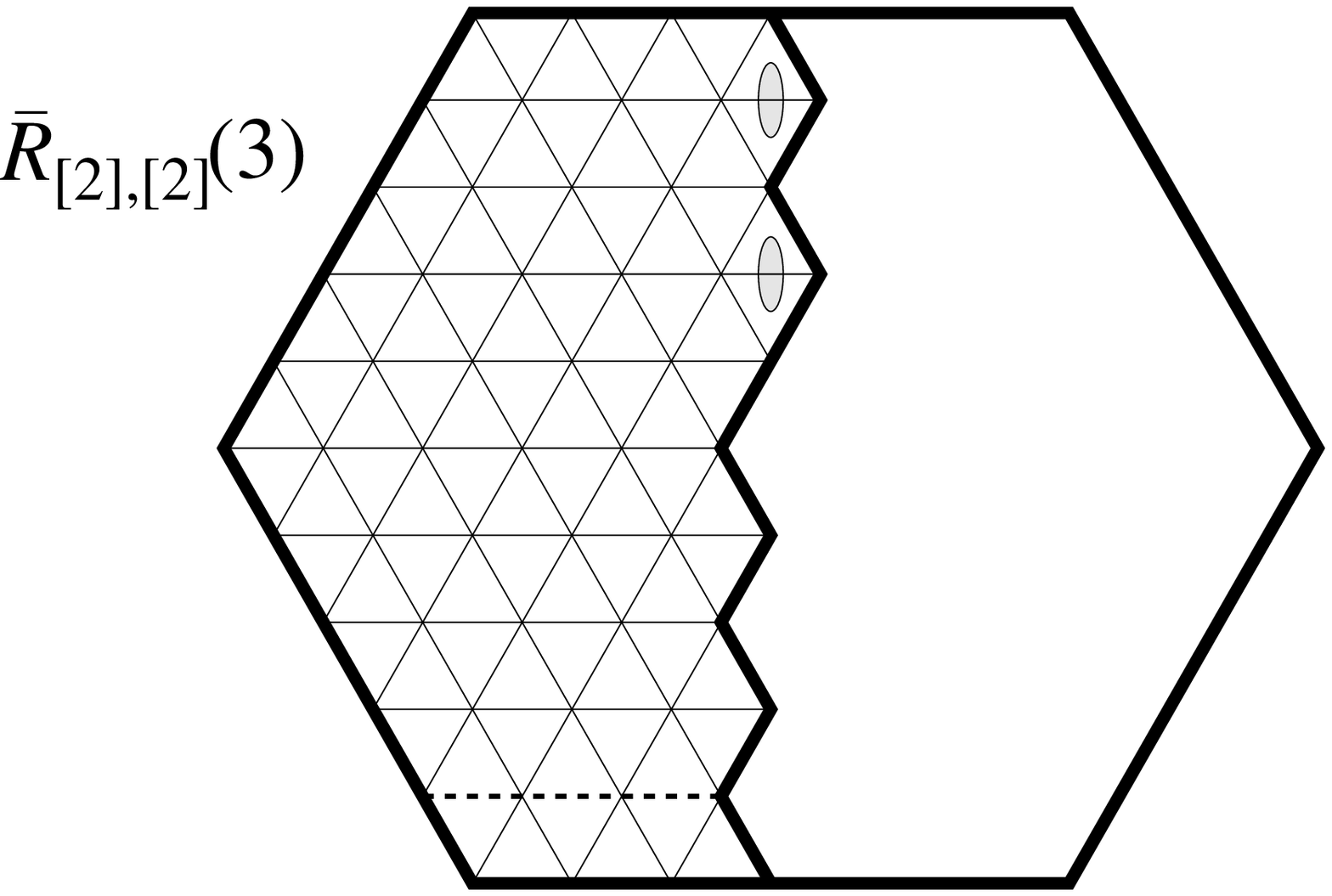}}
\twoline{Figure~7.2. {\rm $a=6$, $b=5$.}}{Figure~7.3. {\rm $a=5$, $b=6$.}}
\endinsert

$$\align
L(H(2x,2x,2y+1))&=2^xL(R_{[x],[x-1]}(y))\tag7.5\\
L(H(2x+1,2x+1,2y))&=2^xL(\bar{R}_{[x],[x]}(y))\tag7.6$$
\endalign$$
(see Figures 7.2 and 7.3).

By Proposition 2.1 of \cite{8}, (7.4)--(7.6) provide product formulas for
$SC(a,a,b)$, and these are easily seen to agree with (7.1)--(7.3). $\square$

\bigskip
{\bf Acknowledgments.} The first author would like to thank John Stembridge, 
David Robbins and Dennis Stanton for stimulating conversations and for their interest 
in this work.

\mysec{References}
{\openup 1\jot \frenchspacing\raggedbottom
\roster

\myref{1}
  T. Amdeberhan, Lewis strikes again!, electronic manuscript dated 1997 (available
at http://www.math.temple.edu/$\sim$tewodros/programs/kradet.ps).
\myref{2}
  G. E. Andrews, Plane partitions, III: The weak Macdonald conjecture, 
{\it Invent. Math.} {\bf 53} (1979), 193--225.  
\myref{3}
  G. E. Andrews, Plane Partitions, V: The T.S.S.C.P.P. conjecture, {\it J.
Comb. Theory Ser. A} {\bf 66} (1994), 28--39.
\myref{4}
  G. E. Andrews and W. H. Burge, Determinant identities, {\it Pacific J. Math.} {\bf
158} (1993), 1--14.
\myref{5}
  G. E. Andrews and D. W. Stanton, Determinants in plane partitions enumeration,
preprint (available at http://www.math.umn.edu/~stanton/PAPERS/det.ps).
\myref{6}
  M. Ciucu, Enumeration of perfect matchings in graphs with reflective symmetry, 
{\it J. Comb. Theory Ser. A} {\bf 77} (1997), 67--97.
\myref{7}
  M. Ciucu, The equivalence between enumerating cyclically symmetric, 
self-com\-ple\-men\-ta\-ry and totally symmetric, self-complementary plane partitions, 
{\it J. Comb. Theory Ser. A}, to appear (available
at http://www.math.ias.edu/$\sim$ciucu/csts.ps).
\myref{8}
  M. Ciucu, Plane partitions I: a generalization of MacMahon's
formula, preprint 
(available at 
http://www.math.ias.edu/$\sim$ciucu/genmac.ps).
\myref{9}
  G. David and C. Tomei, The problem of the calissons, {\it Amer. Math. 
Monthly} {\bf 96} (1989), 429--431.
\myref{10}
  I. M. Gessel and X. Viennot, Binomial determinants, paths, and hook length formulae,
{\it Adv. in Math.} {\bf 58} (1985), 300--321.
\myref{11}
  C. G. J. Jacobi, De formatione et proprietatibus determinantium, {\it J. Reine
Angew. Math.} {\bf 22} (1841), 285--318.
\myref{12}
  C. Krattenthaler, Determinant identities and a generalization of the number of 
totally symmetric self-complementary plane partitions, {\it Elect. J. Combin.} 
{\bf 4} No. 1 (1997), R27.
\myref{13}
  G. Kuperberg, Symmetries of plane partitions and the permanent-de\-ter\-mi\-nant
method, {\it J. Comb. Theory Ser. A} {\bf 68} (1994), 115--151.
\myref{14}
  P. A. MacMahon, Memoir on the theory of the partition of numbers---Part V. Partitions
in two-dimensional space, {\it Phil. Trans. R. S.}, 1911, A.
\myref{15}
  W. H. Mills, D. P. Robbins and H. Rumsey Jr., Proof of the Macdonald conjecture, 
{\it Invent. Math.} {\bf 66} (1982), 73--87. 
\myref{16}
  W. H. Mills, D. P. Robbins and H. Rumsey Jr., Enumeration of a symmetry class of
plane partitions, {\it Discrete Math.} {\bf 67} (1987), 43--55. 
\myref{17}
  R. A. Proctor, Odd symplectic groups, {\it Invent. Math.} {\bf 92} (1988),
307--332.
\myref{18}
  D. P. Robbins, The story of $1$, $2$, $7$, $42$, $429$, $7436,\dotsc,$ {\it Math.
Intelligencer} {\bf 13} (1991), No. 2, 12--19.
\myref{19}
  R. P. Stanley, Symmetries of plane partitions, {\it J. Comb. Theory Ser. A} 
{\bf 43} (1986), 103--113.
\myref{20}
  J. R. Stembridge, Nonintersecting paths, Pfaffians and plane partitions, 
{\it Adv. in Math.} {\bf 83} (1990), 96--131.
\myref{21}
  J. R. Stembridge, The enumeration of totally symmetric plane partitions,
{\it Adv. in Math.} {\bf 111} (1995), 227--243.

\endroster\par}

\enddocument